\documentclass[10pt,twoside]{article}
\usepackage{amsmath}
\usepackage{amsfonts}
\usepackage{amsmath}
\usepackage{amssymb}
\usepackage{fancyhdr}
\usepackage{latexsym}
\usepackage{bbding}
\usepackage{mathrsfs}
\usepackage{wasysym}

\usepackage{multicol,graphics}

\newcommand{\diva}{{\rm div}}
\newcommand{\bff}{{{\bf f}}}

\newtheorem{theorem}{Theorem}[section]
\newtheorem{lemma}[theorem]{Lemma}

\numberwithin{equation}{section}

\allowdisplaybreaks

 \makeatletter
\begin{document}
\title{Global Well-posedness of the  Stochastic Kuramoto- Sivashinsky
Equation with Multiplicative Noise\thanks{Research supported by
the National Natural Science Foundation of China (No. 10771223) and
NSF grant 1025422.} \thanks{E-mail addresses: wuwei5430@yahoo.cn
(Wu); cuisb@yahoo.com.cn (Cui); duan@iit.edu (Duan)}}
\author{{\small   Wei Wu$^{1,2}$}
 \quad {\small  Shangbin Cui$^1$}
 \quad {\small Jinqiao Duan$^3$}
\\
{\small $^{1}$ Department of Mathematics, Sun Yat-Sen University,
Guangzhou, Guandong 510275,   China} \\
{\small $^{2}$ Department of Mathematics, Linyi Normal University,
Linyi, Shandong 276005,   China}\\
{\small $^{3}$ Department of Applied Mathematics, Illinois Institute
of Technology, Chicago, IL 60616, USA}}
\date{}
\maketitle

\begin{abstract}
Global well-posedness of the initial-boundary value problem for the
stochastic Kuramoto-Sivashinsky equation in a bounded domain $D$
with a multiplicative noise is studied. It is shown that under
suitable sufficient conditions, for any initial data $u_0\in
L^2(D\times \Omega)$ this problem has a unique global solution $u$
in the space $L^2(\Omega,C([0,T],L^2({D})))$ for any $T>0$, and
the solution map $u_0\mapsto u$ is Lipschitz continuous.

\textbf{Keywords}: Kuramoto-Sivashinsky equation; stochastic partial differential equations;
multiplicative noise; well-posedness.

\end{abstract}

\section{Introduction}
\noindent The deterministic Kuramoto-Sivashinsky equation was
independently proposed by Kuramoto \cite{Kur} and Sivashinsky
\cite{Siv} as a model describing the instability and turbulence of
wave fronts in chemical reaction and the laminar flames. It also has
many applications in other fields of physics, chemistry, and
biology. We refer the reader to see \cite{Bia}, \cite{Guo},
\cite{Ios}, \cite{Tad}, \cite{Tem}, \cite{Zhang} and references
cited therein for the study on the deterministic
Kuramoto-Sivashinsky equation and its generalizations.

In the present paper, we consider the following initial-boundary value problem of the
stochastic generalized Kuramoto-Sivashinsky equation driven by a
multiplicative noise:
\begin{equation}
\label{1.1} \left\{
\begin{array}{l}
\partial_t u+\Delta^2 u+\Delta u+\diva\bff(u)=\dot{V}_t,\quad x\in{D},\quad t>0,\\
u|_{\partial{D}}=\Delta u|_{\partial{D}}=0,\quad t>0,\\
u|_{t=0}=u_0,\quad x\in{D}.
\end{array}\right.
\end{equation}
Here ${D}$ is a bounded domain in $\mathbb{R}^d$ with a smooth
boundary, $\bff$ is a given $d$-vector function, $\dot{V}_t$ is a
multiplicative noise (see (1.2) below), and $u_0$ is a given initial
$L^2(D)$-valued random variable. The noise is defined in a probability space $(\Omega,\mathscr{F},P)$ and usually we omit the   dependence on
samples $\omega \in \Omega$ in various variables.

The stochastic generalized
Kuramoto-Sivashinsky equation (1.1) is a natural extension of the
deterministic Kuramoto-Sivashinsky equation subject to random influences.
In \cite{DuanE} Duan and Ervin considered the stochastic
Kuramoto-Sivashinsky equation with an additive noise. They proved
global well-posedness of the one-dimensional stochastic
Kuramoto-Sivashinsky equation with an additive noise in the $L^2$
space. The purpose of the present paper is to study the stochastic
generalized Kuramoto-Sivashinsky  equation in the case of
multiplicative noise. By using the truncation method combined with
the $L^2$ conservation law of the Kuramoto-Sivashinsky equation, we
shall prove that the stochastic generalized Kuramoto-Sivashinsky
equation (1.1) is globally well-posed in the $L^2$ space. Here we
remark that in \cite{DuanE}, in order to establish $L^2$
well-posedness of the stochastic Kuramoto-Sivashinsky equation with
additive noise, the authors used a variable transformation to
transform the stochastic equation into a deterministic equation with
the sample point variable as a parameter. This method clearly does
not work for the present multiplicative noise case.

During the past twenty years, great advancement has been made to
the study of stochastic partial differential equations. Some
general theories for such equations have been well-established,
see, for instance, \cite{Chow1}, \cite{DPrat}, \cite{DuanW} and
the references cited therein. We note that despite that the
Kuramoto-Sivashinsky equation is a parabolic equation, the general
theory of parabolic stochastic partial differential equations
developed in the above-mentioned literatures does not apply to the
stochastic Kuramoto-Sivashinsky equation with a multiplicative
noise. This is because that in such a theory the nonlinearity is
required to be of the asymptotically linear type, whereas the
Kuramoto-Sivashinsky equation has a quadratic nonlinear term. Note
that for the stochastic nonlinear wave equations with certain
polynomial nonlinear terms, this difficulty can be overcome with
the aid of the $\dot{H}^1$ conservation law of the nonlinear wave
equations (cf. \cite{Chow2}--\cite{Chow4}). For the
Kuramoto-Sivashinsky equation and its generalized forms, we have
only the $L^2$ conservation law but not any other higher-order
conservation laws. It follows that growth condition in the
generalized Kuramoto-Sivashinsky equations is more restrictive
than the nonlinear wave equations. Similar features are possessed
by the stochastic Burgers equations (cf. \cite{Gy1}, \cite{Gy2}
and \cite{TwaZab}) and the stochastic Navier-Stokes equations (cf.
\cite{CapGat}, \cite{CapPes} and \cite{MikRoz}). However, for
these equations, since they are of the second-order, in order to
get well-posedness of the initial and initial-boundary value
problems, the space dimension $d$ is required to be not greater
than $2$, and for the case $d\geq 3$ we have only existence of
weak solutions but not any well-posedness result (cf.
\cite{CapGat}, \cite{CapPes} and \cite{MikRoz}). For the
stochastic Kuramoto-Sivashinsky equation, as we shall see below,
since it is a fourth-order parabolic equation, well-posedness can
be ensured for $d\leq 5$. For discussions on other fourth-order
parabolic equations, such as the stochastic Cahn-Hilliard
equation, we refer the reader to see \cite{CWeb} and the
references cited therein. We also refer the reader to see
\cite{Yang1} and \cite{Yang2} for the study of long-term behavior
of solutions of the stochastic Kuramoto-Sivashinsky equation (with
additive noise).

We make the following assumption on the nonlinearity $\bff$:

{\bf Assumption $(A)$} $\bff(0)=0$, and there exist constants $C>0$
and $p\geq 1$ such that
\begin{align*}
|\bff(u)-\bff(v)|&\leq C(1+|u|+|v|)^{p-1}|u-v|\quad \mbox{for}\;\;
u,v\in\mathbb{R}.
\end{align*}

As for the noise term $\dot{V}_t$, we assume that it has the
following expression:
\begin{align}
  \dot{V}_t=\sigma(t,x,u,\partial  u,\partial^2
  u)\dot{W_t},
\end{align}
where $\sigma$ is a given function, and $W_t$ is a $L^2(D)$ valued
Wiener process (see Section 2.2 for details; top dots denotes the
derivatives in $t$).

We impose the following assumption on noise intensity $\sigma$:

{\bf Assumption $(B)$} There exist constants $C>0$ and
$\varepsilon>0$ such that
\begin{align*}
|\sigma(t,x,u_1,\xi_1,\zeta_1)-\sigma(t,x,u_2,\xi_2,\zeta_2)|&\leq
C(|u_1-u_2|+|\xi_1-\xi_2|)+\varepsilon|\zeta_1-\zeta_2|)\quad
\mbox{for}\;\; u,v\in\mathbb{R}.
\end{align*}
for $t\geq 0$, $x\in D$, $u,v\in\mathbb{R}$,
$\xi_1,\xi_2\in\mathbb{R}^d$ and
$\zeta_1,\zeta_2\in\mathbb{R}^{\frac{1}{2}d(d-1)}$.

Let $R$ be the covariance operator of the Wiener process $W$, and
$r(x,y)$ be its kernel (see Section 2.2 for details). We need the
following assumption on the kernel function:

{\bf Assumption $(C)$} The kernel function $r$ is in $L^\infty(D\times D)$, so that there
exists a constant $C>0$ such that
\begin{align*}
  r(x,y)\leq C \quad \mbox{for} \;\; x,y\in D.
\end{align*}

Let us now present the main result of this paper. We first consider
the special case that $\sigma$ does not depend on derivatives of
$u$. In this case we have the following result:

\begin{theorem}
Let the assumptions $(A)$, $(B)$ and $(C)$ be satisfied. Suppose
further that $\sigma(t,x,u,\xi,\zeta)=\sigma(t,x,u)$ and $1\leq
p\leq 2$ for $1\leq d\leq 5$ and $1\leq p<1+\frac{6}{d}$ for $d\geq
6$. Then the problem $(1.1)$ is globally well-posed in
$L^2(D\times\Omega)$. More precisely, for any $u_0\in
L^2(D\times\Omega)$ the problem $(1.1)$ has a unique solution $u$
such that for any $T>0$, $u\in L^2(\Omega,C([0,T],L^2({D})))$, and
the solution map $u_0\mapsto u$ is a Lipschitz continuous map from
$L^2(D\times\Omega)$ to $L^2(\Omega,C([0,T],L^2({D})))$.
\hfill$\Box$
\end{theorem}

Next we consider the general case. In this case our result is as
follows:

\begin{theorem}
Let the assumptions $(A)$, $(B)$ and $(C)$ be satisfied. Suppose
further that $1\leq p\leq 2$ for $d=1$ and $1\leq p<1+\frac{2}{d}$
for $d\geq 2$. Then there exists $\varepsilon_0>0$ such that if
$|\varepsilon|\leq \varepsilon_0$, then the problem $(1.1)$ is
globally well-posed in $L^2(D\times\Omega)$. More precisely, for any
$u_0\in L^2(D\times\Omega)$ the problem $(1.1)$ has a unique
solution $u$ such that for any $T>0$, $u\in
L^2(\Omega,C([0,T],L^2({D}))\cap L^2([0,T],H^2({D})))$, and the
solution map $u_0\mapsto u$ is a Lipschitz continuous map from
$L^2(D\times\Omega)$ to $L^2(\Omega,C([0,T],L^2({D}))\cap
L^2([0,T],H^2({D})))$. \hfill$\Box$
\end{theorem}

{\bf Remarks}.\ \ (1)\ \ Note that
$L^2(\Omega,C([0,T],L^2({D})))\hookrightarrow
C([0,T],L^2(D\times\Omega))$ (for deterministic $T>0$). This
justifies the notion of ``well-posedness in $L^2(D\times\Omega)$''.

(2) Throughout this paper, for simplicity we only consider the noise
consisting of a simple term $\sigma\dot{W}$. The theorems hold true
for multiple noise terms $\sum_{i=1}^m\sigma_i\dot{W}_i$ with
independent Wiener fields $W_i$'s, provided that each $\sigma_i$
satisfies the same conditions imposed on $\sigma$.

(3) Our concern in this paper is {\em well-posedness} of the problem
(1.1) in suitable function spaces. Theorems 1.1 and 1.2 give
sufficient conditions to ensure global well-posedness of the problem
(1.1) in the space $L^2(D\times\Omega)$. If one is not concerned
with well-posedness but merely interested in existence of a solution
(so that the solution might not be unique and continuously depend on
the initial data), then these sufficient conditions can be weakened.
We shall discuss this problem in a different paper.

The rest   of this paper is organized as follows. In Section 2 we
present some preliminary materials. In Section 3 we present the proof
of Theorem 1.1, and in Section 4 we present the proof of Theorem 1.2.

\section{Preliminaries}

In this section we give some fundamental estimates for integrals
related to the Green's function $G(x,y,t)$ of the the linear partial
differential equation $\partial_tu+\Delta^2u+\Delta u+cu=0$ (in
${D}$) subject to the boundary value conditions
$u|_{\partial{D}}=\Delta u|_{\partial{D}}=0$. We first consider
deterministic integrals, and next consider stochastic integrals.

\subsection{Estimates for deterministic integrals}

Let $\{\lambda_k\}_{k=1}^\infty$ be the sequence of eigenvalues of
the minus Laplace $-\Delta$ on $D$ subject to the homogeneous
Dirichlet boundary condition, where multiple eigenvalues are counted
in their multiplicities. Let $\{\phi_k\}_{k=1}^\infty$ be the
corresponding sequence of eigenfunctions. We assume that they are
suitably chosen so that they form an  orthonormal basis of
$L^2({D})$. Since $\lim_{k\to\infty}\lambda_k(\lambda_k-1)=\infty$,
there exists $c\geq 0$ such that $\mu_k:=\lambda_k(\lambda_k-1)+c>0$
for all $k\in \mathbb{N}$. Choose a such $c$ and fix it. Let
$$
  G(t,x,y)=\sum_{k=1}^\infty\phi_k(x)\phi_k(y)e^{-\mu_kt},
  \quad x,y\in\overline{{D}},\;\; t>0;
  \quad G(x,y,0)=\delta(x-y).
$$
$G$ is the Green's function of the linear partial differential
equation $\partial_tu+\Delta^2u+\Delta u+cu=0$ (in ${D}$) subject to
the boundary value conditions $u|_{\partial{D}}=\Delta
u|_{\partial{D}}=0$. Note that $\min_{k\geq 1}\mu_k>0$.

\begin{lemma}\label{2.1}\ \ For any $\varphi\in L^2({D})$ and $\alpha
\in\mathbb{Z}_+^d$ we have
\begin{align}
  \|\partial_x^\alpha\!\int_{{D}}\!G(t,x,y)\varphi(y)dy\|_{L^2}\leq
  Ct^{-\frac{|\alpha|}{4}}\|\varphi\|_{L^2}\quad \mbox{for}\;\;
  t>0.
\end{align}
\end{lemma}

{\em Proof}:\ \ For simplicity of the notation we denote
$S(t)\varphi(x)=\int_{{D}}\!G(x,y,t)\varphi(y)dy$. We first consider
the case $\alpha=0$. For $\varphi\in L^2({D})$, let
$\varphi=\sum_{k=1}^\infty a_k\phi_k$. Then
$S(t)\varphi=\sum_{k=1}^\infty a_ke^{-\mu_kt}\phi_k$, so that
\begin{align}
  \|S(t)\varphi\|_{L^2}^2
  =\sum_{k=1}^\infty a_k^2e^{-2\mu_kt}\leq e^{-2c_0t}\sum_{k=1}^\infty a_k^2
  = e^{-2c_0t}\|\varphi\|_{L^2}^2
  \quad (\,c_0=\min_{k\geq 1}\mu_k>0),
\end{align}
by which the assertion for the case $\alpha=0$ follows. Next, since
$\Delta S(t)\varphi=-\sum_{k=1}^\infty a_k\lambda_k
e^{-\mu_kt}\phi_k$, we have
$$
  \|\Delta S(t)\varphi\|_{L^2}^2
  =\sum_{k=1}^\infty a_k^2\lambda_k^2e^{-2\mu_kt}.
$$
Since $\mu_k>0$ for all $k$ and
$\lim_{k\to\infty}\frac{\lambda_k^2}{\mu_k}=1$, there exists a
positive constant $C$, actually $C=\max_{k\geq
1}\frac{\lambda_k^2}{\mu_k}$, such that $\lambda_k^2\leq C\mu_k$ for
all $k$. Hence
\begin{align}
  \|\Delta S(t)\varphi\|_{L^2}^2\leq
  C\sum_{k=1}^\infty a_k^2\mu_k e^{-2\mu_kt}\leq
  Ct^{-1}\sum_{k=1}^\infty a_k^2=Ct^{-1}\|\varphi\|_{L^2}^2.
\end{align}
In getting the second last relation we used the elementary
inequality $xe^{-x}\leq 1/e$ (for $x>0$). Since $S(t)\varphi\in
H^2({D})\cap H^1_0({D})$ for all $t>0$, by making use of the
well-known Agmon-Douglis-Nirenberg inequality (in $L^2$ case) and
the estimates (2.2) and (2.3) we have
$$
  \sum_{|\alpha|=2}\|\partial^\alpha S(t)\varphi\|_{L^2}\leq
  C(\|\Delta S(t)\varphi\|_{L^2}+\|S(t)\varphi\|_{L^2})\leq
  Ct^{-\frac{1}{2}}\|\varphi\|_{L^2}.
$$
This proves the assertion for the case $|\alpha|=2$. The case
$|\alpha|=1$ then follows from interpolation. For general
$\alpha\in\mathbb{Z}_+^d$ the proof is similar. We omit the details.
$\quad\Box$
\medskip

\begin{lemma}\label{2.2}\ \ For any $1\leq q\leq 2$,
$\alpha\in\mathbb{Z}_+^d$ and $\varphi\in L^q({D})$ we have
\begin{align}
  \|\partial_x^\alpha\int_{{D}}G(t,x,y)\varphi(y)dy\|_{L^2}\leq
  Ct^{-\frac{d}{4}(\frac{1}{q}-\frac{1}{2})-\frac{|\alpha|}{4}}
  \|\varphi\|_{L^q},
\end{align}
where $C$ is a positive constant depending only on ${D}$, $d$, $q$
and $\alpha$.
\end{lemma}

{\em Proof}:\ \ We only need to give the proof for the case $q=1$,
because the case $q=2$ is ensured by Lemma 2.1, and the rest cases
follow from these two special cases by interpolation. Moreover, we
may assume that $\varphi\in L^1(D)\cap L^2(D)$, because if this is
proved then for general $\varphi\in  L^1(D)$ the desired assertion
then follows from the density of $L^1(D)\cap L^2(D)$ in $L^1(D)$.
For $\varphi\in L^1(D)\cap L^2(D)$ we let $\varphi=\sum_{k=1}^\infty
a_k\phi_k$. Then $S(t)\varphi=\sum_{k=1}^\infty
a_ke^{-\mu_kt}\phi_k$, so that
$\|S(t)\varphi\|_{L^2}^2=\sum_{k=1}^\infty a_k^2e^{-2\mu_kt}$. We
have
$$
  |a_k|=|\int_D\varphi(x)\phi_k(x)dx|\leq
  \|\varphi\|_{L^1}\|\phi_k\|_{L^\infty}\leq
  C\lambda_k^{\frac{d}{4}}\|\varphi\|_{L^1}.
$$
Here we used the inequality $\|\phi_k\|_{L^\infty}\leq
C\lambda_k^{\frac{d}{4}}$, whose simple proof is as follows: Choose
an integer $l$ sufficiently large such that $2l>d/2$. Then by making
use of the Gagliardo-Nirenberg inequality and the equation
$\Delta^l\phi_k=(-\lambda_k)^l\phi_k$ we have
$$
  \|\phi_k\|_{L^\infty}\leq C\|\phi_k\|_{L^2}^{1-\frac{d}{4l}}
  \|\Delta^l\phi_k\|_{L^2}^{\frac{d}{4l}}=C\lambda_k^{\frac{d}{4}}.
$$
Hence
$$
  \|S(t)\varphi\|_{L^2}\leq
  C\|\varphi\|_{L^1}\Big(\sum_{k=1}^\infty\lambda_k^{\frac{d}{2}}e^{-
  2\mu_kt}\Big)^{\frac{1}{2}}\leq Ct^{-\frac{d}{8}}\|\varphi\|_{L^1}.
$$
By a similar argument we see that for any positive integer $l$,
$$
  \|\Delta^lS(t)\varphi\|_{L^2}\leq
  Ct^{-\frac{d}{8}-\frac{l}{2}}\|\varphi\|_{L^1}.
$$
By using again the Gagliardo-Nirenberg inequality we see that for
any $\alpha\in\mathbb{Z}_+^d$,
$$
  \|\partial_x^\alpha S(t)\varphi\|_{L^2}\leq
  Ct^{-\frac{d}{8}-\frac{|\alpha|}{4}}\|\varphi\|_{L^1}.
$$
This proves the desired assertion. $\quad\Box$
\medskip

We shall also need the following preliminary result which follows
from the energy identity for the equation
$\partial_tu+\Delta^2u+\Delta u+cu=f$:

\begin{lemma}\label{2.3} For $\varphi\in L^2({D})$ we have
\begin{align}
  \int_0^t\!\|\!\int_{{D}}\!
  G(t,x,y)\varphi(y)dy\|_{H^2}^2dt
  \leq C\|\varphi\|_{L^2}^2\quad \mbox{for}\;\;
  t>0.
\end{align}
\end{lemma}

{\em Proof}:\ \ We first assume that $\varphi\in H^2({D})\cap
H^1_0({D})$. In this case $u(t):=S(t)\varphi\in C^\infty([0,\infty),
H^2({D})\cap H^1_0({D}))$, so that the following calculations make
sense. By multiplying both sides of the equation
$\partial_tu+\Delta^2u+\Delta u+cu=0$ with $u$ and integrating over
$D$, we see that
$$
  \frac{1}{2}\frac{d}{dt}\|u(t)\|_{L^2}^2
  +\|\Delta u(t)\|_{L^2}^2-\|\nabla u(t)\|_{L^2}^2+c\|u(t)\|_{L^2}^2=0.
$$
It follows that
\begin{align}
  \|u(t)\|_{L^2}^2
  +2\int_0^t[\|\Delta u(\tau)\|_{L^2}^2-\|\nabla u(\tau)\|_{L^2}^2]d\tau
  \leq\|\varphi\|_{L^2}^2.
\end{align}
Since $u(t)\in H^2({D})\cap H^1_0({D})$ for all $t>0$, we have
$$
  \|\nabla u(t)\|_{L^2}^2\leq\|\Delta u(t)\|_{L^2}\|u(t)\|_{L^2}
  \leq\frac{1}{2}\|\Delta u(t)\|_{L^2}^2+
  \frac{1}{2}\|u(t)\|_{L^2}^2.
$$
Hence, from (2.6) we get
$$
  \int_0^t\|\Delta u(\tau)\|_{L^2}^2d\tau
  \leq\|\varphi\|_{L^2}^2+\int_0^t\|u(\tau)\|_{L^2}^2d\tau.
$$
It follows by the Agmon-Douglis-Nirenberg inequality and (2.4) that
\begin{align*}
  \int_0^t\|u(\tau)\|_{H^2}^2d\tau
  \leq & C\int_0^t[\|\Delta
  u(\tau)\|_{L^2}^2d\tau+\|u(\tau)\|_{L^2}^2]d\tau
  \leq C\|\varphi\|_{L^2}^2+C\int_0^t
  e^{-c_0\tau}\|\varphi\|_{L^2}^2d\tau
  \leq C\|\varphi\|_{L^2}^2.
\end{align*}
For general $\varphi\in L^2({D})$ we use approximation. $\quad\Box$
\medskip

\subsection{Estimates for stochastic integrals}

Let $W_t=W_t(x,\omega)$ ($t\geq 0$) be a $L^2(D)$ valued Wiener
process on a probability space $(\Omega,\mathscr{F}P)$, i.e., there
exists a complete normalized orthogonal basis $\{e_k\}_{k=1}^\infty$
of $L^2(D)$, a sequence of positive numbers $\{c_k\}_{k=1}^\infty$
satisfying $\sum_{k=1}^\infty c_k^2<\infty$, and a sequence of
independent, identically distributed standard Brownian motions
$w_t^k=w_t^k(\omega)$ ($k=1,2,\cdots$) on $(\Omega,\mathscr{F},P)$,
such that
$$
  W_t(x,\omega)=\sum_{k=1}^\infty c_k w_t^k(\omega)e_k(x).
$$
By convention, later on we will omit the sample point variable
$\omega$ and simply write $W_t(x,\omega)$ and $w_t^k(\omega)$
respectively as $W_t(x)$ and $w_t^k$. Note that $W_t$ is the limit
of the finite dimensional Wiener process $W_t^n=\sum_{k=1}^n c_k
w_t^ke_k$ in $L^2(\Omega,C([0,T], L^2(D)))$ (for any $T>0$, cf.
\cite{DPrat}), so that it belongs to $L^2(\Omega,C([0,T],L^2(D)))$
(for any $T>0$). Let
$$
  r(x,y)=\sum_{k=1}^\infty c_k^2 e_k(x)e_k(y).
$$
Then $\int\!\!\int_{D\times D}|r(x,y)|^2dxdy=\sum_{k=1}^\infty
c_k^4\leq(\sum_{k=1}^\infty c_k^2)^2<\infty$, i.e., $r\in
L^2(D\times D)$. Moreover, it is clear that $r(x,y)=r(y,x)$, and
$\int\!\!\int_{D\times D}r(x,y)\varphi(x)\varphi(y)dxdy\geq 0$ for
any $\varphi\in L^2(D)$. Hence $r(x,y)$ defines a positive
semi-definite self-adjoint Hilbert-Schmidt operator $R$ on $L^2(D)$:
$$
  (R\varphi)(x)=\int_{D}r(x,y)\varphi(y)dy \quad \mbox{for}\;\;
  \varphi\in L^2(D).
$$
In fact, $R$ is a self-adjoint trace class operator on $L^2(D)$,
with
$$
  \|R\|_{\mathcal{L}_1}={\rm Tr}R=\int_{D}r(x,x)dx=\sum_{k=1}^\infty c_k^2<\infty.
$$
A simple computation shows that
\begin{align*}
EW_t(x)=0 \quad \mbox{and} \quad E\{W_t(x)W_s(y)\}=(t\wedge s)
r(x,y),
\end{align*}
where $t\wedge s=\min\{t,s\}$. Note also that
$r(x,x)=\sum_{k=1}^\infty c_k^2 e_k^2(x)\geq 0$ for $x\in D$.

Let $\{\mathscr{F}_t\}_{t\geq 0}$ be a filtration of the sub
$\sigma$-fields of $\mathscr{F}$, and $u_t=u(t,x,\omega)$ be a
continuous $L^2(D)$-valued $\mathscr{F}_t$-adapted random field
satisfying the condition
\begin{align}
  E\int_0^T\!\|u_t\|_{R}^2dt<\infty, \quad \mbox{where}\;\;
  \|u_t\|_{R}^2=\int_Dr(x,x)|u(t,x,\omega)|^2dx.
\end{align}
Again, by convention later on we omit the sample point variable
$\omega$ in $u(t,x,\omega)$ and simply write it as $u_t=u(t,x)$.

\begin{lemma}\label{2.4} Assume that the condition $(2.7)$ is satisfied. Then
we have the following estimate:
\begin{align}
  E\Big(\sup_{0\leq t\leq T}\|\int_0^t\int_{D}
  G(t-s,x,y)u(y,s)dydW_s(y)\|_{L^2}^2\Big)\leq
  CE\Big(\int_0^T\|u_t\|^2_R dt\Big).
\end{align}
\end{lemma}

{\em Proof}:\ \ This is a corollary of Theorem 6.10 of \cite{DPrat}.
$\quad\Box$
\medskip

\begin{lemma}\label{2.5} Assume that the condition $(2.7)$ is satisfied. Then
we have the following estimate:
\begin{align}
  \sum_{|\alpha|=2}E\Big(\int_0^T\!\|\partial_x^\alpha\!\!\int_0^t\!\!
  \int_{D}G(t-s,x,y)u(s,y)dydW_s(y)\|_{L^2}^2dt\Big)
  \leq CE\Big(\int_0^T\|u\|^2_Rdt\Big).
\end{align}
\end{lemma}

{\em Proof}:\ \ Since
$G(t,x,y)=\sum_{k=1}^\infty\phi_k(x)\phi_k(y)e^{-\mu_kt}$, we have
\begin{align*}
  \Delta\int_0^t\!\!\int_{D}\!G(t-s,x,y)u(s,y)dydW_s(y)=&
  -\sum_{k=1}^\infty \lambda_ke^{-\mu_kt}\phi_k(x)
  \Big(\int_0^t\!\!\int_{D}\!e^{\mu_ks}\phi_k(y)u(s,y)dydW_s(y)\Big),
\end{align*}
so that
\begin{align}
   &E\int_0^T\!\|\Delta\int_0^t\!\!\int_{D}\!G(t-s,x,y)u(s,y)dydW_s(y)\|_{L^2}^2dt
\nonumber\\
  =&E\int_0^T\!\Big\{\sum_{k=1}^\infty \lambda_k^2e^{-2\mu_kt}
  \Big(\int_0^t\!\!\int_{D}\!e^{\mu_ks}\phi_k(y)u(s,y)dydW_s(y)\Big)^2\Big\}dt
\nonumber\\
  =&\sum_{k=1}^\infty \lambda_k^2\int_0^T\!e^{-2\mu_kt}
  E\Big(\int_0^t\!\!\int_{D}\!e^{\mu_ks}\phi_k(y)u(s,y)dydW_s(y)\Big)^2dt
\nonumber\\
  =&\sum_{k=1}^\infty \lambda_k^2\int_0^T\!e^{-2\mu_kt}
  \Big(E\int_0^t\!e^{2\mu_ks}\big(R\phi_ku_s,\phi_ku_s\big)ds\Big)dt
\nonumber\\
  =&\sum_{k=1}^\infty \lambda_k^2
  E\int_0^T\!\Big(\int_s^T\!e^{-2\mu_kt}dt\Big)e^{2\mu_ks}\big(R\phi_ku_s,\phi_ku_s\big)ds
\nonumber\\
  \leq &\frac{1}{2}\sum_{k=1}^\infty\frac{\lambda_k^2}{\mu_k}
  E\int_0^T\!\big(R\phi_ku_s,\phi_ku_s\big)ds
  \leq CE\int_0^T\!\sum_{k=1}^\infty\big(R\phi_ku_s,\phi_ku_s\big)ds
\nonumber\\
  = &CE\int_0^T\!\!\int_{D}\!r(x,x)u^2(s,x)dxds
  \quad (\mbox{because}\;\;\sum_{k=1}^\infty\phi_k(x)\phi_k(y)=\delta(x-y))
\nonumber\\
  = &CE\int_0^T\!\|u(t,\cdot)\|_{R}^2dt.
\end{align}
In getting the third equality we used the following generalized
It\^{o} isometry:
\begin{align*}
  E\Big(\int_0^t\!\!\int_{D}\!u(s,y)dydW_s(y)\Big)^2=&
  E\int_0^t\!\big(Ru(s,\cdot),u(s,\cdot)\big)ds,
\end{align*}
whose proof is an easy exercise of the stochastic integrals. Indeed,
by letting $J_t(x)=\int_0^t\!u(s,y)dW_s(y)$, we have, by the
stochastic Fubini theorem (see \cite{DPrat}), that
\begin{align*}
  E\Big(\int_0^t\!\!\int_{D}\!u(s,y)dydW_s(y)\Big)^2=&
  E\Big(\int_{D}\!\!\int_0^t\!u(s,y)dW_s(y)dy\Big)^2\\
  =&E\{(J_t,1)(J_t,1)\}=E\int_0^t\!\big(Ru(s,\cdot),u(s,\cdot)\big)ds.
\end{align*}
Having proved (2.10), (2.9) follows immediately from the
Agmon-Douglis-Nirenberg inequality. $\quad\Box$

\section{The proof of Theorem 1.1}

In this section we give the proof of Theorem 1.1. We shall use the
truncation method to prove this theorem.

For every integer $N>0$ we consider a truncated problem as follows:
First we choose a mollifier $\eta_{N}$, i.e., $\eta_{N}:
[0,\infty)\longrightarrow [0,1]$ is a $C^{\infty}$ function
satisfying the condition
\begin{align*}
\eta_{N}(r)=\left\{
\begin{array}{l}
1\quad \mbox{for}\;\; 0\leq r\leq N,\\
0\quad \mbox{for}\;\; r\geq 2N.
\end{array}\right.
\end{align*}
For $u\in L^2$, let $S_N u=\eta_N(\|u\|_{L^2})u$ and
$\bff_N(u)=\bff(S_N u)$. The truncated problem takes the form:
\begin{equation}
\label{1.1} \left\{
\begin{array}{l}
\partial_t u+\Delta^2 u+\Delta u+\diva\bff_N(u)=
\sigma(t,x,u)\dot{W}_{t},
\quad x\in{D},\quad t>0,\\
u|_{\partial{D}}=\Delta u|_{\partial{D}}=0,\quad t>0,\\
u|_{t=0}=u_0,\quad x\in{D}.
\end{array}\right.
\end{equation}
Using the Green's function and the Duhamel's formula, we can convert
the above problem into the following equivalent stochastic integral
equation:
\begin{align}
  u(t,x)&=\int_{D} G(t,x,y)u_0(y)dy+c\int_0^t\!\!\int_{D}
  G(t-s,x,y)u(s,y)dyds\nonumber\\
  &\quad+\int_0^t\!\!\int_{D} \nabla
  G(t-s,x,y)\cdot\bff_N(u(s,y))dy ds\nonumber\\
  &\quad+\int_0^t\!\!\int_{D}\!
  G(t\!-\!s,x,y)\sigma(s,y,u(s,y))dydW_s(y).
\end{align}
In what follows we use the Banach fixed point theorem to prove that
the above problem is globally well-posed in $L^2(D\times\Omega)$.

For any $T>0$, let $X_T$ be the set of $L^2({D})$-valued
$\mathscr{F}_t$-adapted continuous random processes $u$ on $[0,T]$
such that the norm
\begin{align*}
  \|u\|_{X_T}=\Big(E\sup_{0\leq t\leq
  T}\|u\|^2_{L^2}\Big)^{\frac{1}{2}}
\end{align*}
is finite, i.e., $X_T$ is the set of $\mathscr{F}_t$-adapted random
processes belonging to $L^2(\Omega,C([0,T],L^2({D}))$. It is evident
that $(X_T,\|\cdot\|_{X_T})$ is a Banach space. For $u\in X_T$, let
$\Gamma u$ be the right-hand side of (3.2). In what follows we prove
that for any $u\in X_T$, $\Gamma u$ is well-defined and belongs to
$X_T$ as well, and the operator $\Gamma:  X_T\to  X_T$ defined in
this way is a contraction mapping provided $T$ is sufficiently
small.

We first note that the assumption $(A)$ ensures that for any $u,v\in
L^2({D})$,
\begin{align}
  \|\bff(u)\|_{L^\frac{2}{p}}\leq & C(1+\|u\|^{p}_{L^2}),
\\
  \|\bff(u)-\bff(v)\|_{L^\frac{2}{p}}\leq &
  C(1+\|u\|_{L^{2}}+\|v\|_{L^{2}})^{p-1}\|u-v\|_{L^{2}}.
\end{align}
Indeed, by using the assumption $(A)$ we have
\begin{align*}
  \|\bff(u)-\bff(v)\|^\frac{2}{p}_{L^\frac{2}{p}}
  &\leq C\int_{D}(1+|u|+|v|)^{\frac{2}{p}(p-1)}|u-v|^{\frac{2}{p}}dx\nonumber\\
  &\leq C\|u-v\|^\frac{2}{p}_{L^2}(1+\|u\|_{L^2}+\|v\|_{L^2})^{\frac{2}{p}(p-1)},
\end{align*}
by which (3.3) and (3.4) immediately follow. We also note that the
assumptions $(B)$ and $(C)$ ensure that there exists some constant
$C>0$ such that for any $u,v\in L^2({D})$,
\begin{align}
  \|\sigma(t,x,u)\|^2_{R}\leq &
  C(1+\|u\|^2_{L^2}),
\\
  \|\sigma(t,x,u)-\sigma(t,x,v) &\|^2_{R}\leq
  C\|u-v\|^2_{L^2}.
\end{align}
Indeed, the assumptions $(C)$ implies that $\|u\|_{R}\leq
C\|u\|_{L^2}$. Hence, by using the assumptions $(B)$ we immediately
obtain these estimates.

By using Lemma 2.1 with $\alpha=0$ we have
\begin{align}
  E\Big(\sup_{0\leq t\leq T}\|\int_{D} G(t,x,y)u_0(y)dy\|^2_{L^2}\Big)\leq C
  E\Big(\|u_0\|^2_{L^2}\Big),
\end{align}
and
\begin{align}
  E\Big(\sup_{0\leq t\leq T}\|\int_0^t\!\!\int_{D} G(t,x,y)u(s,y)dyds\|^2_{L^2}\Big)
  \leq & CE\Big(\int_0^T\!\|u(s,\cdot)\|^2_{L^2}ds\Big)
  \leq CTE\Big(\sup_{0\leq t\leq T}\|u\|^2_{L^2}\Big).
\end{align}
Next, note that by $(3.3)$ we have
\begin{align}
 \|\bff_N(u)\|_{L^\frac{2}{p}}=\|\bff(S_N u)\|_{L^\frac{2}{p}}
 \leq C(1+\|S_Nu\|^p_{L^2})\leq C(N).
\end{align}
Hence, by using Lemma 2.2 with $|\alpha|=1$ and $q=\frac{2}{p}$, and
noticing the fact that the conditions on $p$ ensures that $1\leq
q=\frac{2}{p}\leq2$ and
$\frac{1}{4}\leq\frac{d}{8}(p-1)+\frac{1}{4}<1$, we have
\begin{align*}
  &\|\int_0^t\!\!\int_{D} \nabla G(t-s,x,y)\cdot\bff_N (u(s,y))dy
  ds\|_{L^2}\\
  \leq &
  \int_0^t\|\int_{D}\nabla G(t-s,x,y)\cdot\bff_N (u(s,y))dy\|_{L^2}ds\\
  \leq & C\int_0^t (t-s)^{-\frac{d}{8}(p-1)-\frac{1}{4}}\|\bff_N
  (u(s,y))\|_{L^\frac{2}{p}}ds\\
  \leq & C(N)\int_0^t(t-s)^{-\frac{d}{8}(p-1)-\frac{1}{4}}ds\\
  =& C(N)t^{\frac{3}{4}-\frac{d}{8}(p-1)},
\end{align*}
which yields
\begin{align}
  E\Big(\sup_{0\leq t\leq T}\|\int_0^t\int_{D}\nabla
  G(t-s,x,y)\cdot\bff_N(u(s,y))dyds\|^2_{L^2}\Big)
  \leq
  C(N)T^{\frac{3}{2}-\frac{d}{4}(p-1)}.
\end{align}
For the stochastic integral, by using Lemma 2.4 and (3.5) we have
\begin{align}
  &E\Big(\sup_{0\leq t\leq T}\|\int_0^t\!\!\int_{D}
  G(t-s,x,y)\sigma(u(s,y))dydW_s(y)\|^2_{L^2}\Big)\nonumber\\
  \leq &CE\Big(\int_0^T\|\sigma(u)\|^2_Rdt\Big)
  \leq CE\Big(\int_0^T(1+\|u\|^2_{L^2})dt\Big)\nonumber\\
  \leq &CT\{1+E(\sup_{0\leq t\leq T}\|u\|^2_{L^2})\}.
\end{align}
Combining the inequalities (3.7), (3.8), (3.10) and (3.11), we see
that there exists constant $C(N,T)>0$ such that
\begin{align*}
  \|\Gamma u\|^2_{X_T}\leq C(N,T)\{1+E(\|u_0\|^2_{L^2})+\|u\|^2_{X_T}\}.
\end{align*}
Therefore, the operator $\Gamma$ is well-defined and maps $X_T$ into
itself.

Next, from (3.2) we see that for $u,v\in X_T$,
\begin{align*}
 \Gamma u-\Gamma v & = c\int_0^t\!\!\int_{D}
  G(t-s,x,y)[u(s,y)-v(s,y)]dyds\nonumber\\
  &\quad+\int_0^t\!\!\int_{D} \nabla
  G(t-s,x,y)\cdot[\bff_N(u(s,y))-\bff_N(v(s,y))]dyds\nonumber\\
  &\quad+\int_0^t\!\!\int_{D}G(t-s,x,y)[\sigma(u(s,y))-\sigma(v(s,y))]dydW_s(y).
\end{align*}
By making use of (3.8) we have
\begin{align}
  E\Big(\sup_{0\leq t\leq T}\|\int_0^t\!\!\int_{D}
  G(t,x,y)[u(s,y)-v(s,y)]dyds\|^2_{L^2}\Big)
  \leq & CTE\Big(\sup_{0\leq t\leq T}\|u-v\|^2_{L^2}\Big).
\end{align}
From (3.4) we see that
\begin{align}
  \|\bff_N(u)-\bff_N(v)\|_{L^\frac{2}{p}}&
  =\|\bff(S_N u)-\bff(S_N v)\|_{L^\frac{2}{p}}
\nonumber\\
  &\leq C(1+\|S_Nu\|_{L^{2}}+\|S_Nv\|_{L^{2}})^{p-1}\|S_Nu-S_Nv\|_{L^{2}}
\nonumber\\
  &\leq C(N)\| u- v\|_{L^{2}}.
\end{align}
Using this inequality and a similar argument as in the proof of
(3.9) we get
\begin{align}
  &E\Big(\sup_{0\leq t\leq T}\|\int_0^t\int_{D}\nabla
  G(t-s,x,y)\cdot[\bff_N (u(s,y))-\bff_N (v(s,y))]dyds\|^2_{L^2}\Big)
\nonumber\\
  \leq &CE\Big\{\sup_{0\leq t\leq T}\Big(\int_0^t
  (t-s)^{-\frac{d}{8}(p-1)-\frac{1}{4}}
  \|\bff_N(u)-\bff_N(v)\|_{L^\frac{2}{p}}ds\Big)^2\Big\}
\nonumber\\
  \leq & C(N)E\Big(\sup_{0\leq t\leq T}\| u- v\|_{L^{2}}^2\Big)
  \Big(\int_0^t(t-s)^{-\frac{d}{8}(p-1)-\frac{1}{4}}ds\Big)^2
\nonumber\\
  \leq &C(N)T^{\frac{3}{2}-\frac{d}{4}(p-1)}
  E\Big(\sup_{0\leq t\leq T}\| u- v\|_{L^{2}}^2\Big).
\end{align}
Finally, by Lemma 2.4 and (3.6) we have
\begin{align}
  &E\Big(\sup_{0\leq t\leq T}\|\int_0^t\int_{D}
  G(t-s,x,y)[\sigma(u(s,y))-\sigma(v(s,y))]dyds\|^2_{L^2}\Big)
\nonumber\\
  \leq & CE\Big(\int_0^T\|\sigma(u)-\sigma(v)\|^2_Rdt\Big)
  \leq CE\Big(\int_0^T\|u-v\|^2_{L^2}dt\Big)
\nonumber\\
  \leq & CTE\Big(\sup_{0\leq t\leq T}\|u-v\|^2_{L^2}\Big).
\end{align}
 Combining (3.12), (3.14) and (3.15), we get
\begin{align*}
  \|\Gamma u-\Gamma v\|^2_{X_T}\leq
  C(N)(T^{\frac{3}{2}-\frac{d}{4}(p-1)}+T)\|u-v\|^2_{X_T}.
\end{align*}
Since $\frac{3}{2}-\frac{d}{4}(p-1)>0$, we see that if $T$ is
sufficiently small so that
$C(N)(T^{\frac{3}{2}-\frac{d}{4}(p-1)}+T)<1$, then the operator
$\Gamma$ is a contraction mapping on $X_T$.

By the Banach fixed point theorem, it follows that if $T$ is so
small that $C(N)(T^{\frac{3}{2}-\frac{d}{4}(p-1)}+T)<1$, then the
equation (3.2) has a unique solution in $X_T$. Since $T$ does not
depend on $u_0$, by a classical argument, the solution can be
extended over all the right-half line $[0,\infty)$, i.e. the
truncated problem (3.1) has a unique global solution $u^N(t,x)$.
Moreover, since this solution is obtained by using the Banach fixed
point theorem, we see that the solution map $u_0\mapsto u^N$ is
Lipschitz continuous from $L^2(\Omega,L^2(D))$ to $X_T$ for any
$T>0$. Hence the problem (3.1) is globally well-posed in
$L^2(\Omega,L^2(D))$.

We now introduce a stopping time $\tau_N$ as follows:
\begin{align*}
  \tau_N=\inf\{t>0:~\|u^N(t,\cdot)\|_{L^2}>N\}
\end{align*}
if the set on the right-hand side is nonempty, and set $\tau_N=T$
otherwise. Then, for $t<\tau_N$, $u(t,x)=u^N(t,x)$ is the solution
of the problem (1.1). Since $\tau_N$ is increasing in $N$, we can
define $\tau_\infty=\lim_{N\rightarrow \infty}\tau_N$. For
$t<\tau_\infty$, we have $t<\tau_N$ for some $N>0$, and we define
$u(t,x)=u^N(t,x)$. By uniqueness of the solution of the truncated
problem (3.1), this definition makes sense. Thus we have proved that
there exists a almost everywhere defined function
$\tau_\infty:\Omega\to (0,\infty]$ such that the problem (1.1) has a
solution on $[0,\tau_\infty)\times D$ almost surely in $\Omega$.
This proves local existence of a solution of the problem (1.1).
Moreover, from the above argument we easily see that if
$\tau_\infty<\infty$, then
$$
  \limsup_{t\uparrow\tau_\infty}\|u(t,\cdot)\|_{L^2}=\infty.
$$

For uniqueness, suppose that there is another solution
$\tilde{u}(t,x)$ defined for $t<\tau$ for a stopping time $\tau$,
i.e., $\limsup_{t\uparrow\tau}\|\tilde{u}(t,\cdot)\|_{L^2}=\infty$.
Then $\tau\geq\tau_N$ for any $N>0$, and $\tilde{u}(t,x)=u^N(t,x)$
for $t<\tau_N$, by uniqueness of the solution of the problem (3.1).
It follows that $\tau\geq\tau_\infty$ and $\tilde{u}(t,x)= u(t,x)$
for $t<\tau_\infty$. This further implies that $\tau=\tau_\infty$.
Therefore, the solution of the problem (1.1) is unique.

To obtain a global solution, we only need to prove that for any
finite $T>0$, there exists a corresponding constant $C(T)>0$ such
that
\begin{align}
   E\|u_{T\wedge\tau_N}\|^2_{L^2}\leq C(T).
\end{align}
Here and hereafter we use the notation $u_{t\wedge\tau_N}$ to denote
the value of $u=u^N$ (defined on the time interval $[0,\tau_N)$) at
the time $t\wedge\tau_N$. Indeed, by the Doob's inequality we have
\begin{align*}
E\|u_{T\wedge \tau_N }\|^2_{L^2}\geq E\{I(\tau_N\leq
T)\|u_{T\wedge \tau_N }\|^2_{L^2}\}\geq N^2P\{\tau_N\leq T\},
\end{align*}
where $I$ denotes the indicate function. If $(3.15)$ holds, then we
get
\begin{align*}
P\{\tau_N\leq T\}\leq\frac{C(T)}{N^2}.
\end{align*}
By the Borel-Cantelli lemma, we have
\begin{align*}
P\{\tau_\infty\leq T\}=0,
\end{align*}
and, therefore, $P\{\tau_\infty> T\}=1$ for any $T>0$. Hence
$u(t,x)=\lim_{N\rightarrow \infty}u^N(t,x)$ is a global solution to
the problem (1.1) as claimed. Therefore, it suffices to prove
(3.16).

Since $u_{t\wedge\tau_N}$ is the solution of the problem (3.1) in
the time interval $[0,T\wedge\tau_N)$, by noticing the fact that
$\bff_N(u)=\bff(u)$ in this time interval and using the It\^{o}'s
formula, we get the following equation:
\begin{align*}
\label{3.17} \|u_{t\wedge \tau_N}\|^2_{L^2}=&\|u_0\|^2_{L^2}-2
\int_0^{t\wedge \tau_N}\!\!(\triangle^2 u_s+\triangle u_s+\diva
\bff(u_s),u_s)ds\\
&+2\int_0^{t\wedge\tau_N}\!\!(\sigma(u_s) dW_s,u_s)+
\int_0^{t\wedge\tau_N}\!\!\|\sigma(u_s)\|^2_R ds.\nonumber
\end{align*}
By integral by parts, we have
\begin{align*}
\|u_{t\wedge\tau_N}\|^2_{L^2}=\|u_0\|^2_{L^2}\!\!-2\!\int_0^{t\wedge\tau_N}
  \!\!(\|\triangle u_s\|^2_{L^2}\!-\!\|\nabla u_s\|^2_{L^2})ds+
  2\!\int_0^{t\wedge\tau_N}\!\!\!(\sigma(u_s) dW_s,u_s)+\!
  \int_0^{t\wedge\tau_N}\!\!\!\|\sigma(u_s)\|^2_R
  ds.
\end{align*}
Taking the expectation and using (3.5), we get
\begin{align*}
    E\|u_{t\wedge\tau_N}\|^2_{L^2}=&E\|u_0\|^2_{L^2}-
  2E\int_0^{t\wedge\tau_N}\!\!
  (\|\triangle u_s\|^2_{L^2}-\|\nabla u_s\|^2_{L^2})ds+
  E\int_0^{t\wedge\tau_N}\!\!\|\sigma(u_s)\|^2_R ds\\
  \leq &E\|u_0\|^2_{L^2}-2E
  \int_0^{t\wedge\tau_N}\!\!(\|\triangle u_s\|^2_{L^2}-\|\nabla
  u_s\|^2_{L^2})ds + CE\int_0^{t\wedge
  \tau_N}\!\!(1+\|u_s\|^2_{L^2})ds.
\end{align*}
Since $\|\triangle u_s\|^2_{L^2}-\|\nabla
u_s\|^2_{L^2}\geq\lambda_1(\lambda_1-1)\| u_s\|^2_{L^2}$, we have
\begin{align*}
  E\|u_{t\wedge \tau_N}\|^2_{L^2} &\leq
  E\|u_0\|^2_{L^2}-2\lambda_1(\lambda_1-1) E \int_0^{t\wedge
  \tau_N}\| u_s\|^2_{L^2}ds+CE\int_0^{t\wedge
  \tau_N}(1+\|u_s\|^2_{L^2})ds\nonumber\\
  &\leq CT+E\|u_0\|^2_{L^2}+
  (C+2\lambda_1-2\lambda_1^2)\int_0^{t}E\|u_{s\wedge
  \tau_N}\|^2_{L^2}ds.
\end{align*}
By the Gronwall's lemma, this yields the following estimate:
\begin{align*}
  E\|u_{t\wedge \tau_N}\|^2_{L^2}\leq
  (E\|u_0\|^2_{L^2}+CT)e^{(C+2\lambda_1-2\lambda_1^2)t}\leq C(T),
\end{align*}
where $C(T)$ is a positive constant independent of $N$. Letting
$t=T$, we see that (3.16) follows. This completes the proof of
Theorem 1.1.

\section{The proof of Theorem 1.2}

In this section we give the proof of Theorem 1.2. Again, we shall
use the truncation method to prove this theorem, but we have to use
a different work space.

For every integer $N>0$ let $\bff_N$ be as before. We consider the
following truncated problem:
\begin{equation}
\label{1.1} \left\{
\begin{array}{l}
\partial_t u+\Delta^2 u+\Delta u+\diva\bff_N(u)=
\sigma(t,x,u,\partial_x u,\partial_x^2 u)\dot{W}_{t},
\quad x\in{D},\quad t>0,\\
u|_{\partial{D}}=\Delta u|_{\partial{D}}=0,\quad t>0,\\
u|_{t=0}=u_0,\quad x\in{D}.
\end{array}\right.
\end{equation}
As before, we can convert the above problem into the following
equivalent stochastic integral equation:
\begin{align}
  u(t,x)&=\int_{D} G(t,x,y)u_0(y)dy+c\int_0^t\!\!\int_{D}
  G(t-s,x,y)u(s,y)dyds\nonumber\\
  &\quad+\int_0^t\!\!\int_{D} \nabla
  G(t-s,x,y)\cdot\bff_N(u(s,y))dy ds\nonumber\\
  &\quad+\int_0^t\!\!\int_{D}\!
  G(t\!-\!s,x,y)\sigma(s,y,u(s,y),\partial_y u(s,y),\partial_y^2 u(s,y))dydW_s(y).
\end{align}
In what follows we use the Banach fixed point theorem to prove that
the above problem is globally well-posed in $L^2(D\times\Omega)$.

For any $T>0$, let $Y_T$ be the set of $L^2({D})$-valued
$\mathscr{F}_t$-adapted continuous random processes $u$ on $[0,T]$
such that the norm
\begin{align*}
  \|u\|_{Y_T}=\Big(E\sup_{0\leq t\leq T}\|u\|^2_{L^2}+
  E\!\int_0^T\!\|u\|^2_{H^2}dt\Big)^{\frac{1}{2}}
\end{align*}
is finite, i.e., $Y_T$ is the set of $\mathscr{F}_t$-adapted random
processes belonging to $L^2(\Omega,C([0,T],L^2({D}))\cap
L^2([0,T],H^2(D)))$. It is evident that $(Y_T,\|\cdot\|_{Y_T})$ is a
Banach space. For $u\in Y_T$, let $\Gamma u$ be the right-hand side
of (4.2). In what follows we prove that for any $u\in Y_T$, $\Gamma
u$ is well-defined and belongs to $Y_T$ as well, and the operator
$\Gamma:  Y_T\to Y_T$ defined in this way is a contraction mapping
provided $T$ is sufficiently small.

We first note that the assumptions $(B)$ and $(C)$ ensure that there
exists some constant $C>0$ and $\varepsilon>0$ such that for any
$u,v\in H^2({D})$,
\begin{align}
  \|\sigma(t,x,u,\partial_x u,\partial^2_x u)\|^2_{R}\leq &
  C(1+\|u\|^2_{L^2})+\varepsilon \|u\|^2_{H^2},
\\
  \|\sigma(t,x,u,\partial_x u,\partial^2_x u)-\sigma(t,x,v,\partial_x v,\partial^2_x v) &\|^2_{R}\leq
  C\|u-v\|^2_{L^2}+\varepsilon \|u-v\|^2_{H^2}.
\end{align}
By using Lemma 2.4 and (4.3) we have
\begin{align}
  &E\Big(\sup_{0\leq t\leq T}\|\int_0^t\!\!\int_{D}
  G(t-s,x,y)\sigma(u(s,y),\partial_x u(s,y),\partial^2_x u(s,y))dydW_s(y)\|^2_{L^2}\Big)
\nonumber\\
  \leq &CE\Big(\int_0^T\|\sigma(u, \partial_x u, \partial^2_x u )\|^2_Rdt\Big)
  \leq E\Big(\int_0^TC(1+\|u\|^2_{L^2})+C(\varepsilon)\|u\|^2_{H^2}dt\Big)
\nonumber\\
  \leq &C(T,\varepsilon)(1+E\sup_{0\leq t\leq T}\|u\|^2_{L^2}+
  E\int_0^T\|u\|^2_{H^2}dt)\leq C(T,\varepsilon)(1+\|u\|^2_{Y_T}).
\end{align}
Combing this with the estimates (3.7), (3.8) and (3.10) in Section
3, we see that there exists constant $C(N,T,\varepsilon)>0$ such
that
\begin{align}
  E\sup_{0\leq t\leq T}\|\Gamma u\|^2_{L^2}\leq
  C(N,T,\varepsilon)\{1+E(\|u_0\|^2_{L^2})+\|u\|^2_{Y_T}\}.
\end{align}
Next, by Lemma 2.3 we have
\begin{align}
  E\int_0^t\!\|\!\int_{{D}}\!
  G(x,y,t)u_0(y)dy\|_{H^2}^2dt
  \leq CE\|u_0\|_{L^2}^2.
\end{align}
Moreover, by using Lemma 2.1 with $|\alpha|=0,2$ we have
\begin{align*}
  &\|\int_0^t\int_{D}  G(t-s,x,y)u(s,y)dy ds\|_{H^2}
  \leq C\int_0^t \|\int_{D}  G(t-s,x,y)u(s,y)dy \|_{H^2} ds
\nonumber\\
  \leq& C\int_0^t\{1+(t-s)^{-\frac{1}{2}}\}\|u\|_{L^2}ds
  \leq C(t+t^\frac{1}{2})\sup_{0\leq s\leq t}
  \|u\|_{L^2},
\end{align*}
so that
\begin{align}
\label{}
  E\int_0^T\|\int_0^t\int_{D}  G(t-s,x,y)u(s,y)dyds
  \|^2_{H^2}dt\leq CT^\frac{3}{2}(1+T^\frac{1}{2})
  E\sup_{0\leq t\leq T}\|u\|^2_{L^2}.
\end{align}
Similarly, by using Lemma 2.2 with $|\alpha|=1,3$ and
$q=\frac{2}{p}$, we have
\begin{align*}
  &\|\int_0^t\int_{D}\nabla G(t-s,x,y)\cdot\bff_N(u(s,y))dyds
  \|_{H^2}
\nonumber\\
  \leq & C\int_0^t \|\int_{D}\nabla G(t-s,x,y)\cdot
  \bff_N(u(s,y))dy\|_{H^2}ds
\nonumber\\
  =& C\int_0^t \|\int_{D} (I-\Delta)\nabla G(t-s,x,y)\cdot
  \bff_N (u(s,y))dy \|_{L^2} ds
\nonumber\\
  \leq & C\int_0^t\{(t-s)^{-\frac{d}{8}(p-1)-\frac{1}{4}}+
  (t-s)^{-\frac{d}{8}(p-1)-\frac{3}{4}}\}
  \|\bff_N(u)\|_{L^{\frac{2}{p}}}ds
\nonumber\\
  \leq & C(N)(t^{\frac{3}{4}-\frac{d}{8}(p-1)}+t^{\frac{1}{4}-\frac{d}{8}(p-1)}),
\end{align*}
so that
\begin{align}
\label{}
  E\int_0^T\|\int_0^t\int_{D}\nabla G(t-s,x,y)\cdot
  \bff_N (u(s,y))dyds\|^2_{H^2}dt
  \leq C(N)(1+T)T^{\frac{3}{2}-\frac{d}{4}(p-1)}.
\end{align}
For the stochastic integral, by using Lemma 2.5 and (4.3) we have
\begin{align}
  &E\int_0^T\|\int_0^t\!\!\int_{D}
  G(t-s,x,y)\sigma(u(s,y),\partial_x u(s,y),\partial^2_x u(s,y))dydW_s(y)\|^2_{H^2}dt
\nonumber\\
  \leq &CE\int_0^T\|\sigma(u,\partial_x u,\partial^2_x u)\|^2_Rdt
  \leq E\Big( \int_0^TC(1+\|u\|^2_{L^2})+C(\varepsilon) \|u\|^2_{H^2}dt\Big)
\nonumber\\
  \leq &C(T,\varepsilon)(1+E\sup_{0\leq t\leq T}\|u\|^2_{L^2}+
  E\int_0^T\|u\|^2_{H^2}dt)\leq C(T,\varepsilon)(1+\|u\|^2_{Y_T}).
\end{align}
Combining (4.6)--(4.10), we see that there exists constant
$C(N,T,\varepsilon)>0$ such that
\begin{align*}
  \|\Gamma u\|^2_{Y_T}\leq C(N,T,\varepsilon)\{1+E(\|u_0\|^2_{L^2})+\|u\|^2_{Y_T}\}.
\end{align*}
Therefore, the operator $\Gamma$ is well-defined and maps $Y_T$ into
itself.

Next, from (4.2) we see that for any $u,v\in Y_T$,
\begin{align*}
 \Gamma u-\Gamma v & = c\int_0^t\!\!\int_{D}
  G(t-s,x,y)[u(s,y)-v(s,y)]dyds\nonumber\\
  &\quad+\int_0^t\!\!\int_{D} \nabla
  G(t-s,x,y)\cdot[\bff_N(u(s,y))-\bff_N(v(s,y))]dyds\nonumber\\
  &\quad+\int_0^t\!\!\int_{D}G(t-s,x,y)[\sigma(s,y,u,\partial_y u,\partial_y^2 u)
  -\sigma(s,y,v,\partial_y v,\partial_y^2 v)]dydW_s(y).
\end{align*}
Thus
\begin{align}
  E\sup_{0\leq t\leq T}&\|\Gamma u-\Gamma v\|^2_{L^2}
  =E\sup_{0\leq t\leq T}\| c\int_0^t\!\!\int_{D}
  G(t-s,x,y)[u(s,y)-v(s,y)]dyds\nonumber\\
  &\quad+\int_0^t\!\!\int_{D}\nabla
  G(t-s,x,y)\cdot[\bff_N(u(s,y))-\bff_N(v(s,y))]dyds\nonumber\\
  &\quad+\int_0^t\!\!\int_{D}G(t-s,x,y)[\sigma(s,y,u,\partial_y u,\partial_y^2 u)
  -\sigma(s,y,v,\partial_y v,\partial_y^2 v)]dydW_s(y)\|^2_{L^2}\nonumber\\
  &\leq CE\sup_{0\leq t\leq T}\{\|\int_0^t\!\!\int_{D}
  G(t-s,x,y)[u(s,y)-v(s,y)]dyds\|^2_{L^2}\nonumber\\
  &\quad+\|\int_0^t\!\!\int_{D} \nabla
  G(t-s,x,y)\cdot[\bff_N(u(s,y))-\bff_N(v(s,y))]dyds\|^2_{L^2}\nonumber\\
  &\quad+\|\int_0^t\!\!\int_{D}G(t-s,x,y)[\sigma(s,y,u,\partial_y u,\partial_y^2 u)
  -\sigma(s,y,v,\partial_y v,\partial_y^2 v)]dydW_s(y)\|^2_{L^2}\},
\end{align}
and
\begin{align}
  E\!\int_0^T\! &\|\Gamma u-\Gamma v\|^2_{H^2}dt=E\!\int_0^T\!\|
  c\int_0^t\!\!\int_{D}
  G(t-s,x,y)[u(s,y)-v(s,y)]dyds\nonumber\\
  &\quad+\int_0^t\!\!\int_{D} \nabla
  G(t-s,x,y)\cdot[\bff_N(u(s,y))-\bff_N(v(s,y))]dyds\nonumber\\
  &\quad+\int_0^t\!\!\int_{D}G(t-s,x,y)[\sigma(s,y,u,\partial_y u,\partial_y^2 u)
  -\sigma(s,y,v,\partial_y v,\partial_y^2 v)]dydW_s(y)\|^2_{H^2}dt\nonumber\\
  &\leq CE\!\int_0^T\!\{\|\int_0^t\!\!\int_{D}
  G(t-s,x,y)[u(s,y)-v(s,y)]dyds\|^2_{H^2}\nonumber\\
  &\quad+\|\int_0^t\!\!\int_{D} \nabla
  G(t-s,x,y)\cdot[\bff_N(u(s,y))-\bff_N(v(s,y))]dyds\|^2_{H^2}\nonumber\\
  &\quad+\|\int_0^t\!\!\int_{D}G(t-s,x,y)[\sigma(s,y,u,\partial_y u,\partial_y^2 u)
  -\sigma(s,y,v,\partial_y v,\partial_y^2 v)]dydW_s(y)\|^2_{H^2}\}dt.
\end{align}
By (3.12) we have
\begin{align}
  E\sup_{0\leq t\leq T}\{\|\int_0^t\!\!\int_{D}
  G(t-s,x,y)[u(s,y)-v(s,y)]dyds\|^2_{L^2}\}
  \leq CT E\sup_{0\leq t\leq T}\|u-v\|^2_{L^2},
\end{align}
and by (3.14) we have
\begin{align}
  &E\sup_{0\leq t\leq T}\{\|\int_0^t\!\!\int_{D} \nabla
  G(t-s,x,y)\cdot[\bff_N(u(s,y))-\bff_N(v(s,y))]dyds\|^2_{L^2}\}\nonumber\\
  &\quad\leq C(N)T^{\frac{3}{2}-\frac{d}{4}(p-1)}E\sup_{0\leq t\leq
  T}\|u-v\|^2_{L^2}.
\end{align}
Moreover, by a similar argument as in the proof of (3.15) and but
using (4.4) instead of (3.6) we have
\begin{align}
  &E\sup_{0\leq t\leq
  T}\{\|\int_0^t\!\!\int_{D}G(t-s,x,y)[\sigma(s,y,u,\partial_y u,\partial_y^2 u)
  -\sigma(s,y,v,\partial_y v,\partial_y^2 v)]dydW_s(y)\|^2_{L^2}\}\nonumber\\
  &\quad\leq CTE\sup_{0\leq t\leq
  T}\|u-v\|^2_{L^2}+C\varepsilon
  E\!\int_0^T\!\|u-v\|^2_{H^2}dt,
\end{align}
and by (4.8) we have
\begin{align}
  E\int_0^T\|\int_0^t\int_{D}G(t-s,x,y)[u(s,y)-v(s,y)]dyds
  \|^2_{H^2}dt\leq CT^\frac{3}{2}(T+1)E\sup_{0\leq t\leq T}\|u-v\|^2_{L^2}.
\end{align}
In addition, by a similar argument as in the proof of  (4.9) but
using (3.14) instead of (3.9), we have
\begin{align}
\label{}
 &E\int_0^T\|\int_0^t\int_{D}\nabla G(t-s,x,y)\cdot
 [\bff_N(u(s,y))-\bff_N(u(s,y))]dyds\|^2_{H^2}dt\nonumber\\
 &\leq C(N)(1+T)T^{\frac{3}{2}-\frac{d}{4}(p-1)}
 E\sup_{0\leq t\leq T}\|u-v\|^2_{L^2}.
\end{align}
Finally, by Lemma 2.5 and (4.4) we have
\begin{align}
  &E\int_0^T\|\int_0^t\!\!\int_{D}
  G(t-s,x,y)[\sigma(s,y,u,\partial_y u,\partial_y^2 u)
  -\sigma(s,y,v,\partial_y v,\partial_y^2 v)]dydW_s(y)\|^2_{H^2}dt
\nonumber\\
  &\leq CE\int_0^T\!\|\sigma(u,\partial_x u,\partial^2_x u)-
  \sigma(v,\partial_x v,\partial^2_x v)\|^2_Rdt
  \leq CE\int_0^T\!\Big(C\|u-v\|^2_{L^2}+\varepsilon\|u-v\|^2_{H^2}\Big)dt
\nonumber\\
  &\leq CTE\sup_{0\leq t\leq T}\|u-v\|^2_{L^2}+
  C\varepsilon E\int_0^T\|u-v\|^2_{H^2}dt.
\end{align}
Combing (4.11)--(4.18), we get
\begin{align*}
  \|\Gamma u-\Gamma v\|^2_{Y_T}
  \leq & C(N)\{T+(1+T)T^{\frac{3}{2}-\frac{d}{4}(p-1)}
  +(1+T)T^{\frac{3}{2}}+T
  +T^{\frac{3}{2}-\frac{d}{4}(p-1)}+T\}\nonumber\\
  &\times E\sup_{0\leq t\leq T}\| u-v\|^2_{L^2}+
  C\varepsilon
  E\!\int_0^T\!\|u-v\|^2_{H^2}dt.
\end{align*}
It follows that if $T>0$ and $\varepsilon>0$ are so small that
\begin{align*}
  C(N)\{T+(1+T)T^{\frac{3}{2}-\frac{d}{4}(p-1)}
  +(1+T)T^{\frac{3}{2}}+T+T^{\frac{3}{2}-\frac{d}{4}(p-1)}+T\}<1
  \quad \mbox{and} \quad
  C\varepsilon<1,
\end{align*}
then we have
$$\|\Gamma u-\Gamma v\|_{Y_T}\leq\delta\| u-v\|_{Y_T}$$
for some $\delta\in (0,1)$ depending on $T$ and $\varepsilon$, i.e.,
$\Gamma$ is a contraction mapping in $Y_T$. Hence, by a similar
argument as in the proof of Theorem 1.1, we see that the desired
assertion follows. This completes the proof of Theorem 1.2.


\medskip


\begin{thebibliography}{99}

\bibitem{Bia} H. A. Biagioni, J. L. Bona, R. J. Iorio and M. Scialom,
 On the Korteweg-de Vries-Kuramoto-Sivashinsky equation. \textit{Adv.
 Diff. Equa.}, \textbf{1}(1996), pp. 1--20.
\bibitem{CapGat} M. Capi\'{n}ski and D, Gatarek, Stochastic
 equations in Hilbert spaces and applications to Navier--Stokes equations
 in any dimension, \textit{J. Funct. Anal.}, \textbf{126}(1994), pp. 26--35.
\bibitem{CapPes} M. Capi\'{n}ski and S. Peszatb, On the existence of a
solution to stochastic Navier--Stokes equations, \textit{Nonlinear
Anal.},  \textbf{44}(2001), pp. 141--177.
\bibitem{CWeb} C. Cardon-Weber, Cahn-Hilliard stochastic equation:
 existence of the solution and of its density, \textit{Bernoulli},
 \textbf{7}(2001), no. 5, pp. 777--816.
\bibitem{Chow1} P.-L. Chow, \textit{Stochastic Partial Differential Equations},
 Chapman \& Hall/CRC, \textbf{}(2007).
\bibitem{Chow2} P.-L. Chow, Stochastic wave equations with polynomial
 nonlinearity, \textit{Ann. Appl. Probab.}, \textbf{12}(2002),
 361--381.
\bibitem{Chow3} P.-L. Chow, Asymptotics of solutions to semilinear
 stochastic wave equations, \textit{Ann. Appl. Probab.},
 \textbf{16}(2006), 757--789.
\bibitem{Chow4} P.-L. Chow, Nonlinear stochastic wave equations : Blow-up
 of second moments in $L^2$ norm, \textit{Ann. Appl. Probab.},
 \textbf{19}(2009), 2039¨C2046.
\bibitem{DPrat} G. Da Prato and J. Zabczyk, \textit{Stochastic Equations
 in Infinite Dimensions}, Cambridge University Press,
 \textbf{}(1992).
\bibitem{DuanE} J. Duan and V. J. Ervin, On the stochastic
 Kuramoto-Sivashinsky equation, \textit{Nonlinear Analysis},
 \textbf{44}(2001), pp. 205-216.
\bibitem{DuanW} J. Duan and W. Wang, \textit{Effective Dynamics of Stochastic
 Partial Differential Equations}, Springer, New York, 2011.
\bibitem{Guo} B. L. Guo and Z. J. Jing, On the generalized Kuramoto-Sivashinsky
 type equations with the dispersive effects, \textit{Annals of
 Mathematical Researches}, \textbf{25}(1992), pp. 1--24.
\bibitem{Gy1} I. Gy\"{o}ngy, Existence and uniqueness results for semilinear
 stochastic partial differential equations, \textit{Stoch. Proc. \&
 Appl.}, \textbf{73}(1998), 271--299.
\bibitem{Gy2} I. Gy\"{o}ngy and D. Nualart, On the stochastic Burgers' equation
 in the line, \textit{Ann. Probab.}, \textbf{27}(1999), 782--802.
\bibitem{Ios} A. Iosevich and J. R. Miller, Dispersive effects in a modified
 Kuramoto-Sivashinsky equation, \textit{Comm. Part. Diff. Equa.},
 \textbf{27}(2002), pp. 2413--2448.
\bibitem{Kur} Y. Kuramoto, Instability and turbulence of wave fronts in
 reaction-diffusion systems. \textit{Prog. Theor. Phys.},
 \textbf{63}(1980), pp. 1885--1903.
\bibitem{MikRoz} R. Mikulevicius and B. L. Rozovskii, Global $L^2$-solutions
 of stochastic Navier-Stokes equations, \textit{Ann. Probab.},
\textbf{33}(2005), 137--176.
\bibitem{Siv} G. I. Sivashinsky, Nonlinear analysis of hydrodynatic
 instabily in laminar flames I: Deriva- tion of basic equations,
 \textit{Acta Astronaut}, \textbf{4}(1977), pp. 1177--1206.
\bibitem{Tad} E. Tadmor, The well-posedness of the Kuramoto-Sivashinsky
 equation. \textit{SIAM J. Math. Anal.}, \textbf{17}(1986), pp.
 884--893.
\bibitem{Tem} R. Temam, \textit{Infinite-Dimensional Dynamical Systems in
 Mechanics and Physics}, 2nd edition, New-York: Springer, 1998.
\bibitem{TwaZab} K. Twardowska and J. Zabczyk, A note on stochastic Burgers'
 system of equations, \textit{Stoch. Anal. Appl.}, \textbf{22}(2004),
 pp. 1641--1670.
\bibitem{Yang1} D. Yang, Random attractors for the stochastic Kuramoto-Sivashinsky
  equation, \textit{Stoch. Anal. Appl.}, \textbf{24}(2006), pp. 1285¨C1303.
\bibitem{Yang2} D. Yang, Dynamics for the stochastic nonlocal Kuramoto-Sivashinsky
  equation, \textit{J. Math. Anal. Appl.}, \textbf{330}(2007), pp.
  550¨C570.
\bibitem{Zhang} L. Zhang, Decay of solutions of the multidimensional
 generalized Kuramoto-Sivashinsky System, \textit{IMA J. Appl. Math.},
 \textbf{50}(1993), pp. 29--42.
\end{thebibliography}
\end{document}